\documentclass{amsart}
\usepackage{amsmath,amsthm}
\usepackage{amsfonts,amssymb}

\newtheorem{thm}{Theorem}[section]
\newtheorem{cor}[thm]{Corollary}
\newtheorem{lem}[thm]{Lemma}
\newtheorem{prop}[thm]{Proposition}

\newtheorem{defn}[thm]{Definition}

\theoremstyle{remark}

\def\hb{\hfill$\Box$}
\def\ta{\theta}
\def\al{{\alpha}}
\def\be{{\beta}}
\def\da{{\delta}}
\def\Bl{\Bigl}
\def\Br{\Bigr}
\def\f{\frac}
\def\vi{\varphi}
 \def\a{{\alpha}}
 \def\b{{\beta}}
 \def\g{{\gamma}}
 \def\k{{\kappa}}
 \def\t{{\theta}}
 \def\l{{\lambda}}
 \def\d{{\delta}}
 \def\o{{\omega}}
 \def\s{{\sigma}}
 \def\la{{\langle}}
 \def\ra{{\rangle}}
\def\va{\varepsilon}
 \def\ve{{\varepsilon}}

 \def\Kb{{\mathbf K}}

 \def\CH{{\mathcal H}}

 \def\CO{{\mathcal O}}

 \def\CS{{\mathcal S}}
 
 \def\CV{{\mathcal V}}

 \def\NN{{\mathbb N}}

 \def\RR{{\mathbb R}}
 \def\ZZ{{\mathbb Z}}
        
        \def\proj{\operatorname{proj}}
        \def\supp{\operatorname{supp}}

\def\df{\displaystyle\frac}

\def\be{\beta}
\def\dmin{\displaystyle \min}
\def\sa{\sigma}
\newcommand{\wt}{\widetilde}

\begin{document}

\title[Ces\`aro means of orthogonal expansions]
{Ces\`aro means of orthogonal expansions in several variables}
\author{Feng Dai}
\address{Department of Mathematical and Statistical Sciences\\
University of Alberta\\, Edmonton, Alberta T6G 2G1, Canada.}
\email{dfeng@math.ualberta.ca}
\author{Yuan Xu}
\address{Department of Mathematics\\ University of Oregon\\
    Eugene, Oregon 97403-1222.}\email{yuan@math.uoregon.edu}

\date{\today}
\keywords{Ces\`aro means, $h$-harmonics, sphere, orthogonal polynomials,
ball, simplex}
\subjclass{33C50, 42B08, 42C10}
\thanks{The first  author  was partially supported  by the NSERC Canada
under grant G121211001. The second author was partially supported by the
NSF under Grant DMS-0604056}

\begin{abstract}
Ces\`aro $(C,\delta)$ means are studied for orthogonal expansions with
respect to the weight function $\prod_{i=1}^{d}|x_i|^{2\k_i}$ on the unit
sphere, and for the corresponding weight functions on the unit ball and
the Jacobi weight on the simplex. A sharp pointwise estimate is established
for the $(C,\d)$ kernel with $\d > -1$ and for the kernel of the projection
operator, which allows us to derive the exact order for the norm of the
 Ces\`aro means and the projection operator on these domains.
\end{abstract}

\maketitle

\section{Introduction}
\setcounter{equation}{0}

It is well known that Ces\`aro $(C,\d)$ means of the Jacobi
polynomial expansions with respect to the weight function
$(1-t)^\a (1+t)^\b$ on $[-1,1]$ converges uniformly if and only if
$\delta > \max\{\a,\b\}+ 1/2$ (\cite{Szego}, \cite[p. 78,
Corollary 18.11] {Ch-Mu}). Recently, results as such have been
extended to orthogonal expansions in several variables  (see
\cite{DX, LX, X97} and the references therein). In the present
paper we study orthogonal expansions and their Ces\`aro $(C,\d)$
means with respect to the weight functions
\begin{equation}  \label{weight}
  h_\k(x): = \prod_{i=1}^{d+1} |x_i|^{\k_i}, \qquad \k_i \ge 0,
\end{equation}
on the unit sphere $S^{d} = \{x: \|x\| =1\} \subset \RR^{d+1}$,
where $\|x\|$ denotes the Euclidean norm, as well as similar
problems for orthogonal expansions on the unit ball with respect
to the weight function
\begin{equation} \label{weightB}
 W_\k^B(x): = \prod_{i=1}^{d} |x_i|^{\k_i}(1-\|x\|^2)^{\k_{d+1}-1/2},
      \qquad \k_i \ge 0,
\end{equation}
on the unit ball $B^d = \{x: \|x\| \le 1\} \subset \RR^d$, and for the orthogonal
expansion with respect to the weight function
\begin{equation}  \label{weightT}
 W_\k^T(x) := \prod_{i=1}^{d} x_i ^{\k_i-1/2}(1-|x|)^{\k_{d+1}-1/2},
      \qquad \k_i \ge 0,
\end{equation}
on the simplex $T^d =\{x:x_1\ge 0, \ldots, x_d\ge 0, 1-|x| \ge 0\}$,
where $|x| := x_1+ \cdots + x_d$.

A homogeneous polynomial orthogonal with respect to $h_\k^2$ on
the unit sphere is called an $h$-harmonic. The theory of
$h$-harmonics is developed by Dunkl (see \cite{DX} and the
references therein) for a family of weight functions invariant
under a finite reflection group, of which $h_\k$ in \eqref{weight}
is the simplest example of the group $\ZZ_2^{d+1}$. Let
$\CH_n^{d}(h_\k^2)$ denote the space of spherical $h$-harmonics of
degree $n$. It is known that $\dim
\CH_n^d(h_\k^2)=\binom{n+d+1}{n}-\binom{n+d-1}{n-2}$. The usual
Hilbert space theory shows that
$$
L^2(h_k^2,S^{d}) = \sum_{n=0}^\infty \CH_n^{d}(h_\k^2),  \qquad
       f = \sum_{n=0}^\infty \proj_n(h_\k^2; f),
$$
where $\proj_n(h_\k^2) : L^2(h_\k^2;S^{d}) \mapsto \CH_n^{d}(h_\k^2)$ is the
projection operator, which can be written as an integral operator
\begin{equation} \label{projection}
   \proj_n(h_\k^2; f, x) = a_\k \int_{S^{d}} f(y) P_n(h_\k^2 ;x,y) h_\k^2(y)
   d\o(y),\   \    x\in S^d,
\end{equation}
where $d\o(y)$ denotes the usual Lebesgue measure on $S^d$, and
$P_n(h_\k^2)$ is the reproducing kernel of $\CH_n^{d}(h_\k^2)$.

A fundamental result for our study is the following compact
expression of this kernel (\cite{D2, X97b} or \cite[[p. 202]{DX})
\begin{align} \label{proj-Kernel}
 P_n(h_\k^2; x,y) = c_\k  \frac{n+\l_\k}{\l_\k} \int_{[-1,1]^{d+1}}
      C_n^{\l_\k}(u(x,y,t)) \prod_{i=1}^{d+1}(1+t_i) (1-t_i^2)^{\k_i-1}dt,
\end{align}
where $C_n^\l$ is the Gegenbauer polynomial of degree $n$,
\begin{equation} \label{lambda}
  \l_\k = |\k|+\frac{d-1}{2},  \quad|\k|=\sum_{j=1}^{d+1}\k_j, \quad u(x,y,t)=
      x_1y_1t_1+\ldots + x_{d+1}y_{d+1}t_{d+1},
\end{equation}
and  $c_\k$ is the normalization constant of the weight function
$\prod_{i=1}^d(1+t_i) (1-t_i^2)^{\k_i-1}$.

For $\delta > -1$, the Ces\`aro $(C,\delta)$ means of the $h$-harmonic
expansion is defined by
$$
 S_n^\d (h_\k^2;f,x) : = (A_n^\delta)^{-1}\sum_{k=0}^n 
       A_{n-k}^\delta \proj_n(h_\k^2; f,x),
 \qquad A_{n-k}^\delta = \binom{n-k+\delta}{n-k}.
$$
The case $\d = -1$ can be considered as $\proj_n(h_\k^2; f)$ itself.
Evidently the $(C,\d)$ means can be written as an integral against a kernel,
$K_n^\d(h_\k^2;x,y)$; that is,
$$
S_n^\d (h_\k^2;f,x) : = a_\k \int_{S^{d}} f(y) K_n^\d(h_\k^2 ;x,y) h_\k^2(y) d\o(y),
$$
where $K_n^\d(h_\k^2)$ is the $(C,\d)$ mean of the kernel
$P_n(h_\k^2)$ and $a_\k$ is the normalization constant $a_\k =
1/\int_{S^{d}} h_\k^2d\o$. Many results on $h$-harmonic expansions
have been developed by now. In the following we only state those
results that are essential for our study, refer to \cite{DX} for
the background and refer to \cite{LX} for results on $(C,\d)$
means. Let $P_n^{(\a,\b)}$ denote the $n$-th Jacobi polynomial,
which is the orthogonal polynomial with respect to the weight
function
$$
w^{(\a,\b)}(t) = (1-t)^\a (1+t)^\b, \qquad t \in [-1,1]
$$
with the usual normalization (\cite{Szego}). The Gegenbauer polynomial
$C_n^{\l}$ corresponds to $\a= \be = \l -1/2$, although the
normalization constant is different \cite[p. 80]{Szego}. Let
$K_n^\d(w^{(\a,\b)};s,t)$ denote the $(C,\d)$ means of the kernel of
the Jacobi expansion on $[-1,1]$. Then it follows from
\eqref{proj-Kernel} that
\begin{align} \label{kernel}
K_n^\d(h_\k^2;x,y) = &\ c_\k \int_{[-1,1]^{d+1}} K_n^\d(w^{(\l_\k-\f12,\l_\k-\f12)};1,
          u(x,y,t)) \\
    &\qquad \qquad\qquad \times \prod_{i=1}^{d+1}(1+t_i) (1-t_i^2)^{\k_i-1}dt. \notag
\end{align}
If some $\kappa_i  =0$, then the formula holds under the limit relation
\begin{equation}\label{1-8-07}
 \lim_{\lambda \to 0} c_\lambda \int_{-1}^1 f(t) (1-t)^{\lambda -1} dt
  = \frac12[f(1) + f(-1)] .
\end{equation}

Similar results hold for orthogonal expansions on the unit ball
$B^d$ and on the simplex $T^d$. Let $\Omega^d$ and $W$ denote
either $B^d$ and $W_\k^B$ or $T^d$ and $W_\k^T$, respectively. Let
$\CV_n^d(W)$ denote the space of orthogonal polynomials of degree
$n$ and $\proj_n(W): L^2(W) \mapsto \CV^d_n(W)$ the orthogonal
projection. The Ces\`aro $(C,\d)$ means of the orthogonal
expansion with respect to $W$ are defined as the $(C,\d)$ means of
$\proj_n(W;f)$. These means can also be written as integral
operators,
$$
    S_n^\d(W; f,x) = a_\k^\Omega \int_{\Omega} f(y)  \Kb_n^\d(W; x,y)
       W(y) dy,
$$
where the kernel $\Kb_n^\d(W)$ is the $(C,\d)$ mean of the
reproducing kernels of $\CV^d_n(W)$ and $a_k^\Omega$ is the
normalization constant of $W$ on $\Omega$. There is a close
relation between orthogonal expansions with respect to $W_\k^B$ on
$B^d$ and  the $h$-harmonic expansions with respect to $h_\k^2$ on
$S^d$. In particular, it is known that
\begin{align} \label{kernelB}
\Kb_n^\d (W_{\k}^B; x, y) = & \frac{1}{2}\left [
       K_n^\d(h_{\k}^2;(x,x_{d+1}),(y, y_{d+1})) \right. \\
   & \qquad \left.   +   K_n^\da(h_{\k}^2;(x,x_{d+1}),(y, -y_{d+1}))\right]  \notag
\end{align}
where $x_{d+1}= \sqrt{1-\|x\|^2}$, $y_{d+1}= \sqrt{1-\|y\|^2}$.
Because of this identity, the pointwise estimate of the kernel
$\Kb^\da_n(W_\k^B;x,y)$ can be deduced from that of
$K_n^\da(h_\k^2;x,y)$. There is also a close relation between
orthogonal polynomials on $B^d$ and those on $T^d$, but it is a
relation that involves a transform akin to the quadratic transform
between the Jacobi polynomials and the Gegenbauer polynomials (see
\cite[(4.3.4) and (4.1.5)]{Szego}). The kernel for $W_\k^T$ on
$T^d$ is more complicated as it is given by
\begin{align} \label{kernelT}
\Kb_n^\d (W_\k^T;x,y)  =
  & c_\k \int_{[-1,1]^{d+1}} K_n^\d \left(w^{(|\k| +\frac{d-2}{2},-\frac12)}; 1,
        2z(x,y,t)^2-1\right) \\
   & \qquad \times    \prod_{i=1}^{d+1} (1-t_i^2)^{\k_i-1} d t, \notag
\end{align}
where
$$
z(x,y,t) =\sqrt{x_1 y_1}\, t_1 + \ldots + \sqrt{x_d y_d}\, t_d+ \sqrt{1-|x|}
\sqrt{1-|y|}\, t_{d+1}.
$$

In the case of $d =1$, the weight function $W_\k^T$ becomes the Jacobi
weight $w^{(\k_1 -\frac{1}{2}, \k_1 -\frac{1}{2})}(t)$, so that our results reduce
to the result for Jacboi expansions. The weight function $W_\k^B$ when
$d=1$ becomes the weight function
$$
  w_{\k_2,\k_1} (t) = |t|^{2 \k_1} (1- t^2)^{\k_2-1/2}, \qquad \k_i \ge 0, \quad
   t \in [-1,1],
$$
whose corresponding orthogonal polynomials, $C_n^{(\k_1,\k_2)}$, are
called generalized Gegenbauer polynomials, and they can be expressed
in terms of Jacobi polynomials,
\begin{align}\label{G-Gegen}
\begin{split}
C_{2n}^{(\lambda ,\mu )}(t) &=\frac{\left( \lambda +\mu \right)_{n}}
{\left( \mu +\frac{1}{2}\right)_{n}} P_{n}^{(\lambda -1/2,\mu-1/2)}
(2t^{2}-1), \\
C_{2n+1}^{(\lambda ,\mu )}(t) &=\frac{\left( \lambda +\mu
\right)_{n+1}} {\left( \mu
+\frac{1}{2}\right)_{n+1}}tP_{n}^{(\lambda -1/2,\mu+1/2)}
(2t^{2}-1),
\end{split}
\end{align} where $(a)_n =a(a+1)\cdots (a+n-1)$.
Furthermore, let $\wt C_n^{(\l,\mu)}$ denote the orthonormal generalized
Gegenbauer polynomial; then we have (\cite{X97b})
\begin{align} \label{prodGG}
 & \tilde C_n^{(\lambda,\mu)}(x) \tilde C_n^{(\lambda,\mu)}(y)
    = \frac{n + \lambda + \mu}{\lambda + \mu} c_\lambda c_\mu
        \int_{-1}^1 \int_{-1}^1 \\
  & \quad C_n^{\lambda+\mu}
      (t x y + s\sqrt{1-x^2} \sqrt{1-y^2} )
  (1+t)(1-t^2)^{\mu -1} (1-s^2)^{\lambda -1} dt ds, \notag
\end{align}
which plays an essential role in our proof of various lower bounds.

The convergence of the Ces\`aro means with respect to $h_\k^2$ was
first proved in \cite{X97} under the condition $\delta > |\k|
+\frac{d-1}{2}$. The critical index of the $(C,\d)$ means turns
out to be $\delta > |\k| + \frac{d-1}{2} - \min_{1\le i \le d+1}
\k_i$,  which was proved in \cite{LX} together with similar
results for orthogonal expansions on $B^d$ and on $T^d$.  The main
ingredient of the proof is a sharp pointwise estimate for the
$(C,\d)$ kernel function that was established for $\delta \ge
(d-1)/2$. The derivation of the estimate in \cite{LX} is elaborate
and lengthy, and cannot be extended to $\d < (d-1)/2$. Moreover,
the estimate for the kernel $K_n^\d(W_\k^T;x,y)$ on the simplex
was established under an additional restriction on $\k$, so that
the result on $T^d$ was incomplete.

In the present paper we will establish the pointwise estimate of
the $(C,\d)$ kernel for all $\d > -1$, as well as for the kernel of
the orthogonal projection operator itself, with a much more elegant
proof. As a consequence, we are able to determine the exact order
of the norm of the $(C,\d)$ means for all $\d \ge -1$, including
the projection operator and the partial sum operator, for the
orthogonal expansions on the sphere, the ball, and the simplex.
The deviation of the main estimate on the kernel function
$K_n^\d(h_\k^2;x,y)$ comes down to estimate a multiple integral of
the Jacobi polynomial that has boundary singularities, which in
fact holds for even weaker condition than what is needed for $\d
\ge -1$; both the proof and the result could be useful for other
problems. The sharpness of the norm relies on a lower bound for a
double integral of Jacobi polynomials, which was established in
\cite{LX} in the case of critical index. We will extend this lower
bound to $\d \ge -1$ by using asymptotic expansion of integrals,
which gives a proof that is not only more general but also more
elegant even in the case of critical index.

The paper is organized as follows. The main results are stated and proved
in the following section, assuming the estimates of the kernel. The pointwise
estimate of the kernel is established in Section 3. The lower bound estimate
is established in Section 4.


\section{Main results}
\setcounter{equation}{0}

Throughout this paper we denote by $c$ a generic constant that
may depend on fixed parameters such as $\k$ and $p$, whose value may
change from line to line. Furthermore we write $A \sim B$ if $A \ge c B$
and $B \ge c A$.

\subsection{Orthogonal  expansion on the sphere}
The main estimate of the kernel function is as follows:

\begin{thm} \label{thm:S-estimate}
Let  $x=(x_1,\cdots, x_{d+1})\in S^d$ and $y=(y_1,\cdots, y_{d+1})
\in S^{d}$. Then for  $\d
> -1$,
\begin{align} \label{eq:S-est}
|K_n^\delta(h_\kappa^2; x,y)| \le \;  c & \left[
\frac{ \prod_{j=1}^{d+1}(|x_jy_j|+n^{-1}\|\bar x- \bar y\|+n^{-2})^{-\kappa_j}}
 {n^{\delta -(d-1)/2} (\|\bar x - \bar y\|+ n^{-1})^{\delta+(d+1)/2} }
 \right. \\
 & {\hskip .5in}   +  \left.
\frac{ \prod_{j=1}^{d+1}(|x_jy_j|+ \|\bar x- \bar y\|^2+n^{-2})^{-\kappa_j}}
     {n (\|\bar x - \bar y\|+n^{-1})^{d+1}} \right ],\notag
\end{align} where  $\bar z = (|z_1|, \ldots,  |z_{d+1}|)$ for $z=(z_1,\cdots, z_{d+1})\in S^d$.
Furthermore, for the kernel of projection operator,
\begin{equation} \label{eq:S-estP}
|P_n(h_\kappa^2; x,y)| \le \, c
\frac{ \prod_{j=1}^{d+1}(|x_jy_j|+n^{-1}\|\bar x- \bar y\|+n^{-2})^{-\kappa_j} }
 {n^{-(d-1)/2} (\|\bar x - \bar y\|+ n^{-1})^{(d-1)/2 } }.
\end{equation}
\end{thm}

In the following we take the convention that in the  case $\d =
-1$, $S_n^\d(h_\k^2;f)$ is understood to be just $\proj_n(h_\k^2; f)$.
This pointwise estimate was proved in \cite{LX} for $\d \ge
(d-1)/2$. For $1\le p \le \infty$ let $\|\cdot\|_{\k,p}$ denote
the usual $L^p(h_\k^2; S^{d})$ norm, where in the case of $p =
\infty$ we consider $C(S^d)$, the space of continuous functions
with uniform norm $\|f\|_{\k,\infty} := \|f\|_\infty$.  Let
$\|S_n^\d(h_\k^2)\|_{\k,p}$ denote the operator norm of
$S_n^\d(h_\k^2)$ as an operator from $L^p(h_\k^2;S^d)$ to
$L^p(h_\k^2;S^d)$. As a consequence of the main estimate, we can
prove the following:

\begin{thm} \label{thm:2.2}
Let $\d > -1$ and define
$$
   \s_\k := \tfrac{d-1}{2}+|\kappa|- \min_{1 \le i \le d+1} \kappa_i.
$$
Then for $p = 1$ and $p= \infty$,
$$
    \|\proj_n(h_\k^2)\|_{k,p} \sim n^{\s_\k} \quad \hbox{and}\quad
   \|S_n^\d(h_\k^2)\|_{k,p} \sim  \begin{cases}
          1, & \d > \s_\k  \\
          \log n, & \d = \s_\k\\
           n^{- \delta + \s_\k},  & -1<\d < \s_\k
\end{cases}.
$$
In particular, $S_n^\delta(h_\k^2;f)$ converges in $L^p(h_\k^2;S^d)$ for all
$1 \le p \le \infty$ if and only if $\d > \s_\k$.
\end{thm}

The last statement means that $\s_\k$ is the critical index of the $(C,\d)$ means,
which was proved earlier in \cite{LX}. The results for $\d < \s_\k$ are new. Let
us mention the particular two interesting cases. One is $\d =0$ for which
$S_n^\d$ becomes the partial sum operator
$$
       S_n(h_\k^2; f) = \sum_{j=0}^n \proj_j^\k f,
$$
which is the best approximation to $f$ in $L^2(h_\k^2;S^d)$. The other case
is the projection operator itself.

\begin{cor}
For $p=1$ or $\infty$, $\|S_n(h_\k^2) \|_{\k,p} \sim \|\proj_n(h_\k^2)\|_{\k,p}
 \sim n^{\s_\k }$.
\end{cor}

The proof of Theorem \ref{thm:S-estimate} will be given in Section
3. The estimate of the norm $\|S_n^\d (h_\k^2)\|_{\k,p}$  for
$p=1$ and $p=\infty$   in Theorem $\ref{thm:2.2}$  implies that
the same estimate holds for $ 1 < p < \infty$. For $\d > \s_\k$,
this shows that $\|S_n^\d(h_\k^2)\|_{\k,p}$ is bounded for $1 \le
p \le \infty$. For $\d < \s_\k$, however, the estimate is not
sharp. For example, we know that $\|\proj_n(h_\k^2; f)\|_{\k, 2} \le
\|f\|_{\k,2}$.

While the proof of Theorem $\ref{thm:2.2}$ follows along the same
line as that of \cite[Theorem 2.1]{LX}, which concerns only with the
case of the critical index, it is necessary to provide proofs for several
subtle points, especially for the lower bound. Below we shall present
a self-contained proof. The proof is naturally divided into two parts,
one deals with the upper bound of the norm, the other concerns with
the lower bound of the norm.

\medskip\noindent
{\it Proof of Theorem $\ref{thm:2.2}$ (upper bound).}\  \ We shall
prove the upper bound for the norm of $S_n^\d(h_\k^2)$ with $\da>
-1$. The case of projection operator can be treated similarly.

A standard duality argument shows that $\|S_n^\d(h_\k^2)\|_{\k,1}
= \|S_n^\d(h_\k^2)\|_{\k,\infty}$ so that we only need to consider
the case of $\|\cdot\|_{\k,\infty}$ norm, which is given by
\begin{equation} \label{eq:infty-norm}
    \|S_n^\d (h_\k^2)\|_{\k,\infty} = \sup_{x \in S^d} a_\k \int_{S^d}
          |K_n^\d(h_\k^2;x,y)| h_\k^2(y) d\o(y).
\end{equation}

We claim that
\begin{equation}\label{2-4-07}
|K_n^\da (h_\k^2; x, y)| h_\k^2(y) \leq c n^d (1+n \|\bar x - \bar
y\|)^{-\be(\da)},\   \    \  x,y\in S^d,
\end{equation}
with $\be(\da) =\dmin\{ d+1, \da-\sa_\k+d\}$.  Once the claim
($\ref{2-4-07}$) is proven, then we have
\begin{align*}
  \int_{S^d} |K_n^\d(h_\k^2;x,y)| h_\k^2(y) d\o(y) &\, \leq c \,
          n^d\int_0^{\f\pi2} (1+n\ta)^{-\be(\da)} (\sin \t)^{d-1}\, d\ta \\
      &\,    \sim \begin{cases}
 1, & \d > \s_\k  \\
          \log n, & \d = \s_\k\\
           n^{- \delta + \sigma_\k},  & -1< \d < \s_\k
 \end{cases},
\end{align*}
which together with ($\ref{eq:infty-norm}$) will  give  the
desired upper bound of $\|S_n^\d (h_\k^2)\|_{\k,p}$.

For the proof of ($\ref{2-4-07}$),  we shall use Theorem
\ref{thm:S-estimate}. Without loss of generality we may assume
$|x_1|=\max_{1\leq j\leq d+1}|x_j|$.  Set
$$
I_j(x, y): =(|x_jy_j|+n^{-1}\|\bar x- \bar
y\|+n^{-2})^{-\kappa_j}|y_j|^{2\k_j},\   \  1\leq j\leq d+1.
$$
Since  $|x_1|=\max_{1\leq j\leq d+1}|x_j|\ge  \f 1{\sqrt{d+1}}$, we  have
$$
I_1(x, y) \le |x_1|^{-\k_1}|y_1|^{\k_1}\leq (d+1)^{\f{\k_1}2}.
$$
For $j \ge 2$, if $|x_j|\ge 2\|\bar x-\bar y\|$ then $|y_j|\le  |x_j|+\|\bar
x-\bar y\|\leq \f32|x_j| $,  and hence
$$
I_j (x,y) \leq |x_jy_j|^{-\k_j}|y_j|^{2\k_j} \leq \left(\tfrac32\right)^{\k_j}
  \leq \left(\tfrac32\right)^{\k_j} (1+n\|\bar x-\bar
y\|)^{\k_j},
$$
whereas if $|x_j|< 2\|\bar x-\bar y\|$ then
$|y_j|\leq |x_j|+\|\bar x-\bar y\|\leq 3\|\bar x-\bar y\| $,  and hence
$$
I_j(x,y)\leq (n^{-1} \|\bar x-\bar y\|+n^{-2})^{-\k_j} (3\|\bar x-\bar
y\|)^{2\k_j}\leq 9^{\k_j} (1+n\|\bar x-\bar y\|)^{\k_j}.
$$
Consequently, it follows that
$$
  \prod_{j=1}^{d+1} I_j (x, y) \leq c \prod_{j=2}^{d+1} I_j (x,y)
    \leq c (1+n\|\bar x-\bar y\|)^{|\k|- \k_1}, 
$$
in which $\k_1$ can be replaced by $\min_{1 \le i \le d+1} \k_i$.
Thus, we obtain
\begin{align} \label{2-5-07}
I(x, y) & :=n^{d} (1+n \|\bar{x}-\bar{y}\|)^{-\da -\f{d+1}2}
\prod_{j=1}^{d+1} I_j(x, y)\\
&\leq c n^{d} (1+n
\|\bar{x}-\bar{y}\|)^{-(\da+d-\sa_\k)}.\notag
\end{align}

Similarly, one can show that for  $1\leq j\leq d+1$,
$$
J_j(x, y): =(|x_jy_j|+\|\bar x-\bar y\|^2+n^{-2})^{-\k_j}
|y_j|^{2\k_j} \leq c,
$$
which implies that
\begin{equation}\label{2-6-07}
J(x,y): = \f { \prod_{j=1}^{d+1}
J_j(x,y)}{ n( n^{-1} +\|\bar x-\bar y\|)^{d+1}}
\leq c n^d (1+ n \|\bar x-\bar y\|)^{-d-1}.
\end{equation}

Since Theorem \ref{thm:S-estimate} shows that
$$
  |K_n^\da (h_\k^2; x, y)| h_\k^2(y)\leq c ( I(x, y)+J(x,y)),
$$
the claim ($\ref{2-4-07}$)  follows by ($\ref{2-5-07}$) and ($\ref{2-6-07}$).
\hb

\medskip\noindent
{\it Proof of Theorem $\ref{thm:2.2}$ (lower bound).}\  \  The lower
bound of the norm $\|S_n^\d (h_\k^2)\|_{\k,p}$ follows from the lower
bound in Theorem \ref{thm:2.4} below.  Here we only consider
the case of projection operator.

Let $e_1=(1,0,\ldots,0), \ldots,  e_{d+1} = (0,\ldots,0,1)$ be the
standard basis of $\RR^{d+1}$. By \eqref{proj-Kernel} and \eqref{prodGG},
$$
 P_n(h_\k^2;x,e_j) = \wt C_n^{(\l_k-\k_j, \k_j)}(1)\wt C_n^{(\l_k-\k_j, \k_j)}(x_j)
    = \frac{n+\l_k}{\l_k} C_n^{(\l_k-\k_j, \k_j)}(x_j),
$$
where the second equal sign follows from \cite[p. 27]{DX}. Consequently,
if $\min_{1\le i \le d+1} \k_i= \k_j$ for $1 \le j \le d+1$, then
\begin{align*}
\|\proj^\k_n(h_\k^2)\|_{\k,1} & \ge a_\k \int_{S^d} |P_n(h_\k^2;x,e_j)| h_\k^2(x) d\o(x) \\
 & = \frac{n+\l_k}{\l_k} a_\k \int_{S^d} \left| C_n^{(\l_k-\k_j, \k_j)}(x_j) \right|
  h_\k^2(x) d \o(x) \\
& \ge \frac{n+\l_k}{\l_k}  c \int_{-1}^1
   \left| C_n^{(\s_\k, \k_j)}(x_j) \right| w_{\s_\k,\k_j}(x_j) dx_j.
\end{align*}
Next we write the last integral as twice of the integral over $[0,1]$, as
justified by \eqref{G-Gegen}, and then change variable $2 x_j^2-1 \mapsto t$.
Using \eqref{G-Gegen} we then conclude that
\begin{align*}
 \|\proj^\k_{2n}(h_\k^2)\|_{\k,1} \ge c  n^{\s_\k+\f12} \int_{-1}^1
   \left |P_n^{(\s_\k-\f12,\k_j-\f12)}(t)\right| w^{(\s_\k-\f12,\k_j-\f12)}(t)  dt
   \sim n^{\s_\k},
\end{align*}
where the last step follows from the classical estimate for the
integral of Jacobi polynomials in \cite[(7.34.1)]{Szego}. The case
of $\proj^\k_{2n+1}(h^2_\k)$ is handled similarly. \qed

\subsection{Orthogonal  expansion on the ball}
The pointwise upper bound of the kernel  $\Kb_n^\d(W_\k^B;x,y)$
can be derived from Theorem \ref{thm:S-estimate} using the
identity \eqref{kernelB}. In fact, for our main results on the
norm of $(C,\d)$ means, we can use \eqref{kernelB} directly. For
$1\le p \le \infty$ let $\|\cdot\|_{W_\k^B,p}$ denote the
$L^p(W_\k^B; B^d)$ norm, where in the case of $p = \infty$ we
consider $C(B^d)$ with uniform norm $\|f\|_{W_\k^B,\infty} :=
\|f\|_\infty$.  Let $\|S_n^\d(W_\k^B)\|_{\k,p}$ denote the
operator norm of $S_n^\d(W_\k^B)$ as an operator from $L^p(W_\k^B;
B^d)$ to $L^p(W_\k^B;B^d)$.

\begin{thm} \label{thm:2.4}
Let $\d > -1$ and define $\s_\k := \tfrac{d-1}{2}+|\kappa|-
\min_{1 \le i \le d+1} \kappa_i.$ Then for $p =1$ or $\infty$,
$$
   \|S_n^\d(W_\k^B)\|_{W_k^B,p} \sim  \begin{cases}
          1, & \d > \s_\k  \\
          \log n, & \d = \s_\k\\
           n^{- \delta + \sigma_\k},  & -1<\d < \s_\k
      \end{cases}. \qquad
$$
In particular, $S_n^\delta(W_\k^B;f)$ converges in $L^p(W_\k^B;B^d)$ for all
$1 \le p \le \infty$ if and only if $\d > \s_\k$. Furthermore,
$$
   \|\proj_n(W_\k^B)\|_{W_k^B,p} \sim  n^{\s_\k}
$$
unless $\min_{1\le i \le d+1} \k_i = \k_{d+1}$ and $n$ is odd, in which case
the norm has an upper bound of $c\,n^{\s_\k}$.
\end{thm}

Again the fact that $\s_\k$ is the critical index of the $(C,\d)$ means was proved
earlier in \cite{LX}. The results for $\d < \s_\k$ are new. Let $S_n(W_\k^B; f)$
denote the partial sum operator
$$
       S_n(W_\k^B; f) = \sum_{j=0}^n \proj_n(W_\k^B; f).
$$

\begin{cor}
For $p =1$ or $\infty$, $\|S_n(W_\k^B) \|_{W_\k^B,p}  \sim n^{\s_\k }$.
\end{cor}

Recall that the weight function $W_\k^B$ becomes $w_{\k_2,\k_1}$ in the
case of $d =1$, so that the results of Theorem  \ref{thm:2.4} and its corollary
hold for the generalized Gegenbauer expansions.  Moreover, let
$K_n^\d(w_{\l,\mu}; s,t)$ denote the $(C,\d)$ kernel for the generalized
Gegaubauer expansion with respect to $w_{\l,\mu}$ and define
\begin{equation} \label{eq:Tn-d}
T_n^\d (w_{\l,\mu};t):= \int_{-1}^1
|K_n^\d(w_{\l,\mu};s,t)|w_{\l,\mu}(s)ds;
\end{equation}
then the following proposition plays an essential role in establishing the lower
bound in Theorem \ref{thm:2.4}.

\begin{prop} \label{prop:Tn-lwbd}
Assume $\mu \ge 0$ and $\d \le \l$.
If $\l \ge \mu$ then
$$
  T_n^\d(w_{\l,\mu}; 1),  T_n^\d (w_{\mu,\l};0) \ge
       c  n^{-\d + \l} \begin{cases}
                         \log n, & \text{if $ \d = \l$,}\\
                          1, & \text{if $-1<\d < \l $}.  
                      \end{cases}
$$
\end{prop}

This proposition will be established in  Section 4. Below we use the
proposition to prove Theorem \ref{thm:2.4}.

\medskip\noindent
{\it Proof of Theorem \ref{thm:2.4}. }
The upper bound of the norm in Theorem \ref{thm:2.4} follows easily from
that of Theorem \ref{thm:2.2} as shown in the proof in \cite[p. 286]{LX}.
For the lower bound estimate, the case $\d > -1$ follows essentially the
the proof in \cite{LX}, which is based on the following inequality
(see \eqref{eq:infty-norm}),
$$
\|S_n^\d(W_\k^B)\|_\infty \ge  a_\k^B \int_{B^d}
  |{\bf K}_n^\d(W_\k^B,y, e)|W_\k^B(y)dy : = \Lambda_n(e),
$$
where $e$ is a fixed point in $B^d$. Let $e_1=(1,0,\ldots,0), \ldots,
e_d = (0,\ldots,0,1)$ be the standard basis of $\RR^d$. Following
\cite[p. 287]{LX}, we have
$$
\Lambda_n(e_j) = c T_n^\d(w_{\l_\k-\k_j,\k_j};1), \quad 1 \le j
\le d,\quad \hbox{and} \quad \Lambda_n(0) = c
T_n^\d(w_{\k_{d+1},\l_\k-\k_{d+1}};0),
$$
from which the lower bound of the norm estimate in Theorem \ref{thm:2.4}
follows from Proposition \ref{prop:Tn-lwbd}.

Next we consider the norm of the projection operator. If
$\min_{1\le i \le d+1} \k_i = \k_j$ for $1 \le j \le d$, then by
\eqref{proj-Kernel} and \eqref{prodGG}
$$
 P_n(W_\k^B;x,e_j) = \wt C_n^{(\l_k-\k_j, \k_j)}(1)\wt C_n^{(\l_k-\k_j, \k_j)}(x_j)
    = \frac{n+\l_k}{\l_k} C_n^{(\l_k-\k_j, \k_j)}(x_j),
$$
so that the proof follows exactly as in the case of lower bound of Theorem
\ref{thm:2.4}. We are left with the case of $\min_{1\le i \le d+1} \k_i= \k_{d+1}$.
In this case,  it follows by  the projection operator version of \eqref{kernelB}
and \eqref{prodGG} that
\begin{equation} \label{GG(0)}
  P_n(W_\k^B;x,0) = \wt C_n^{(\k_{d+1}, \s_\k)}(0)\wt C_n^{(\k_{d+1},\s_k)}(|x|).
\end{equation}
Hence, using the structure constants given in \cite[p. 27]{ DX} and
\eqref{G-Gegen}, we obtain that
$$
   P_{2n}(W_\k^B; x,0) = (-1)^n \frac{2n+\l_k}{\l_k}
   \frac{(\l_k)_n}{(\k_{d+1}+\f12)_n} P_n^{(\k_{d+1}-\f12,\s_\k-\f12)}(2\|x\|^2-1).
$$
Using the polar coordinates and then changing variable $2 r^2 -1 \mapsto t$, it
follows that
\begin{align*}
&  \int_{B^d} \left | P_{2n}(W_\k^B;x,0) \right | W_\k^B(x)dx  \\
&   \qquad  \sim n^{\s_\k + \f12}
      \int_{-1}^1\left| P_n^{(\k_{d+1} -\f12,\s_k-\f12)} (t) \right|
          w^{(\k_{d+1} -\f12,\s_k-\f12)}(t) dt  \sim n^{\s_\k}
\end{align*}
again by \cite[(7.34.1)]{Szego}.
\qed

\medskip

Note that by \eqref{GG(0)} and \eqref{prodGG}, $P_{2n+1}(W_\k^B;x,0) \equiv 0$
so that the above method fails when $\k_{d+1} = \min_{1 \le j \le d+1} \k_j$ and
$n$ is odd.

\subsection{Orthogonal  expansion on the simplex}
As mentioned in the introduction, the pointwise estimate of
$\Kb_n^\d(W_\k^T;x,y)$ is more complicated and it does not follow
directly from that of $K_n^\d(h_\k^2;x,y)$.  To state the result,
we introduce the following notation: for $x=(x_1,\cdots,
x_d),y=(y_1,\cdots, y_d) \in T^d$,
$$
\xi := (\sqrt{x_1}, \ldots, \sqrt{x_d}, \sqrt{x_{d+1}}), \qquad
\zeta := (\sqrt{y_1}, \ldots, \sqrt{y_d}, \sqrt{y_{d+1}})
$$ with $x_{d+1}:=1-|x|$ and
$y_{d+1}:=1-|y|$.  Both of these two are points in $S^d$ as $|x| =
x_1+\cdots +x_{d}$ by definition.

\begin{thm} \label{thm:T-estimate}
Let $\d > -1$. For $x ,y \in T^d$,
\begin{align} \label{eq:T-est}
 |\Kb_n^\delta(W_{\kappa}^T; x,y)| \le & \; c \left [
  \frac{ \prod_{j=1}^{d+1} (\sqrt{x_j y_j} + n^{-1} \|\xi- \zeta\|+
    n^{-2})^{- \kappa_j} }
 { n^{\delta - (d-1)/2} (\|\xi - \zeta\|+ n^{-1})^{\delta  +(d+1)/2}} \right. \cr
&\; \qquad \quad + \left. \frac{ \prod_{j=1}^{d+1} (\sqrt{x_j y_j} +\|\xi- \zeta\|^2 +
  n^{-2})^{- \kappa_j} }{ n (\|\xi - \zeta\|+ n^{-1})^{d+1}} \right ].  \notag
\end{align}
Furthermore, for the kernel of the projection operator,
\begin{equation} \label{eq:T-estP}
 |P_n(W_{\kappa}^T; x,y)| \le \, c
  \frac{ \prod_{j=1}^{d+1} (\sqrt{x_j y_j} + n^{-1} \|\xi- \zeta\|+
    n^{-2})^{- \kappa_j} }
 { n^{- (d-1)/2} (\|\xi - \zeta\|+ n^{-1})^{ (d-1)/2}}
\end{equation}
\end{thm}

This estimate was proved in \cite{LX} for $\d \ge (d-2)/2$ and an additional
restriction on $\k$. As in the case of $B^d$ we let $\|\cdot\|_{W_\k^T,p}$
denote the $L^p(W_\k^T; T^d)$ norm and let $\|S_n^\d(W_\k^T)\|_{\k,p}$
denote the operator norm of $S_n^\d(W_\k^T)$ as an operator from
$L^p(W_\k^T; T^d)$ to $L^p(W_\k^T;T^d)$.

\begin{thm} \label{thm:2.7}
Let $\d > -1$ and define $\s_\k := \tfrac{d-1}{2}+|\kappa|-
\min_{1 \le i \le d+1} \kappa_i.$ Then for $p=1$ or $\infty$,
$$
   \|\proj_n(W_\k^T)\|_{W_k^T,p} \sim n^{\s_\k} \quad\hbox{and}\quad
   \|S_n^\d(W_\k^T)\|_{W_k^T,p} \sim  \begin{cases}
          1, & \d > \s_\k  \\
          \log n, & \d = \s_\k\\
           n^{- \delta + \sigma_\k},  & -1< \d < \s_\k
 \end{cases}.
 $$
In particular, $S_n^\delta(W_\k^T;f)$ converges in $L^p(W_\k^T;T^d)$ for all
$1 \le p \le \infty$ if and only if $\d > \s_\k$.
\end{thm}

The fact that $\s_\k$ is the critical index of the $(C,\d)$ means was proved
in \cite{LX} under an additional condition of
$\sum_{i=1}^{d+1} (2\k_i - \lfloor \k_i\rfloor) \ge 1 + \min_{1 \le i \le d=1}\k_i$.
This restriction is now removed. Let $S_n(W_\k^T; f)$ denote the partial
sum operator of the orthogonal expansion.

\begin{cor}
For $p=1$ or $\infty$, $\|S_n(W_\k^T) \|_{W_\k^T,p} \sim
\|\proj_n(W_\k^T) \|_{W_\k^T,p} \sim n^{\s_\k }$.
\end{cor}

\medskip\noindent
{\it Proof of Theorem \ref{thm:2.7}}. \ \
The proof of the upper bound follows from the proof of \cite[Theorem 2.9]{LX},
which reduces the estimate to the one in Theorem \ref{thm:2.2} for all $\d$
and the same reduction holds also for the projection operator. For the lower
bound estimate, we note that (\cite[p. 290]{LX})
\begin{align*}
  \Kb_n^\d(W_\k^T;x, e_j) & =
       K_n^\d \left(w^{(\l_\k -\k_j-\frac12, \k_j -\frac12)};1, 2x_j-1\right),
                        \quad 1 \le j\le d, \\
  \Kb_n^\d(W_\k^T;x, 0) & =
       K_n^\d \left(w^{(\l_\k -\k_j-\frac12, \k_j -\frac12)};1, 1-2|x|\right);
\end{align*}
and the similar formulas hold for projection operator, in which the right hand
side holds with $P_n(w^{(\a,\b)};s,t):= \wt P_n^{(\a,\b)}(s)\wt P_n^{(\a,\b)}(t)$,
where $\wt P_n^{(\a,\b)}(s)$ is the orthonormal polynomial.  Consequently,
as in \cite{LX}, the lower bound estimate reduces to that of Jacobi expansions
at the point $x=1$, for which the relevant results can be deduced easily
from \cite[Chapt. 9]{Szego} (see Lemma \ref{lem:3.6} below).
\qed

\medskip

The results stated above are for the norm of the operators. For the
pointwise convergence, we have the following result.

\begin{thm} \label{thm:2.9}
Let $f$ be continuous on $T^d$. If $\delta > (d-1)/2$, then the $(C, \delta)$
means $S_n^\d (W_\k^T; f)$ converge to $f$ at every point in the interior
of $T^d$ and, furthermore, the convergence is uniform over any compact
set contained in the interior of $T^d$.
\end{thm}

This theorem was proved in \cite{LX} under the condition
$\sum_{j=1}^{d+1} (\k_i - \lfloor \k_i\rfloor) \ge 1$.  The proof
uses a local estimate of the kernel derived from the main estimate
in Theorem \ref{thm:T-estimate}, hence is valid now for all $\d > (d-1)/2$.
Similar pointwise convergences also hold for $S^d$ and $B^d$,
see \cite{LX}.


\section{Pointwise estimates of the kernels}
\setcounter{equation}{0}

The center piece of the pointwise estimate on the $(C,\d)$ kernel is an
estimate of integrals on Jacobi polynomials. This is presented in the
first subsection, from which the estimate of the kernels will be derived
in the subsequent subsections.

\subsection{Main estimate}
The following theorem contains the key ingredient for our pointwise estimate.

\begin{thm} \label{MainEst}
Assume $\k_j>0$, $a_j\neq 0$ and $\vi_j\in C^\infty[-1,1]$ for
$j=1,2,\cdots, m$. Let $|a|: = \sum_{j=1}^m|a_j|\leq 1$.
If $\al\ge \be$,  $\a \ge |\k|- \f12:=\sum_{j=1}^m\k_j -\f12$ and $|x|+|a|\leq 1$,
then
\begin{align} \label{eq:MainEst}
&  \left | \int_{[-1,1]^m} P_n^{(\al,\be)} \biggl(\sum_{j=1}^m a_j t_j +x \biggr)
\prod_{j=1}^m \vi_j(t_j) (1-t_j^2)^{\k_j-1} \, dt \right | \\
& \quad\qquad
\leq c n^{\al-2|\k|} \frac{\prod_{j=1}^m  (|a_j| + n^{-1}\sqrt{1-|a|-|x|}+n^{-2})^{-\k_j}}
  { \left(1+n \sqrt{1- |a|-|x|}\right)^{\a+\f12-|\k|}}. \notag
\end{align}
\end{thm}

It is well known that the Jacobi polynomials satisfy the following
estimate (\cite[(7.32.5) and (4.1.3)]{Szego}).

\begin{lem} \label{lem:3.2}
For an arbitrary real number $\alpha$ and $t \in [0,1]$,
\begin{equation} \label{Est-Jacobi}
|P_n^{(\alpha,\beta)} (t)| \le c n^{-1/2} (1-t+n^{-2})^{-(\alpha+1/2)/2}.
\end{equation}
The estimate on $[-1,0]$ follows from the fact that $P_n^{(\alpha,\beta)}
(t) = (-1)^nP_n^{(\beta, \alpha)} (-t)$.
\end{lem}
The Jacobi polynomials also satisfy the following identity
\begin{equation} \label{D-Jacobi}
    P_n^{(a+\frac12, b+\frac12)}(y) = \frac{2}{n+a+b+1}
      \frac{d}{dy} P_{n+1}^{(a-\frac12, b-\frac12)}(y).
\end{equation}
Hence, in terms of the power of $n$, \eqref{Est-Jacobi} is most
useful for $\a < 1$.  In order to use the inequality effectively,
we give the following definition.

\begin{defn}
Let $n, v\in\mathbb{N}_0$, $\mu,r \in \RR$ with $r >0$. Assume
$|\rho|+r\leq 1$. A function $f:[-r,r]\to \mathbb{R}$ is said
to be in  class $ \mathcal {S}^v_n ( \rho, r, \mu)$, if there exist
functions $F_j$, $j=0,1,\cdots, v$ on $[-r,r]$ such that $F_j^{(j)}(x)=f(x)$,
$x\in [-r,r]$, $0\leq j\leq v$, and
\begin{equation}\label{2-1}
|F_j(x)| \le c n^{-2j} \left(1+n \sqrt{1-|\rho+x|}\right)^{-\mu-\f12+j},\ \ \
x\in [-r,r],\   \    \  j=0,1,\cdots, v.
\end{equation}
\end{defn}

We note that $n^{-\a} P_n^{(\a,\b)} \in \CS_n^v(0,1,\a)$ for all
$v \in \NN_0$ by \eqref{Est-Jacobi} and \eqref{D-Jacobi}.

\begin{lem} \label{lem:3.4}Assume $\d >0$ and $0<|a|\leq r$.
Let $f\in \mathcal {S}^v_n ( \rho,r, \mu)$ with  $v\ge |\mu|+2\da
+\f32$, and let $\xi\in C^\infty [-1,1]$ be such that
$\text{supp}\ \xi \subset [-\f12, 1]$. Then
\begin{align}\label{1-3}
 \left | \int_{-1}^1 f(at+x)(1-t)^{\da-1} \xi(t)\, dt\right |
   \le  c n^{-2\da} |a|^{-\da} (1+n \sqrt{1-A})^{-\mu-\f12+\da} ,
\end{align}
where $A:= |\rho+a+x|$ and $|x|\leq r-|a|$.
\end{lem}

\begin{proof}
To simplify the notation, we define
$$
B:= \frac{1+n\sqrt{ 1-A}}{2n^2|a|}.
$$
First we claim that for $t\in[1-B,1]$,
\begin{equation}\label{2-2}
   1+ n \sqrt{ 1-|at+x+\rho|}\sim 1+n \sqrt{ 1-A}=2n^2|a|B.
\end{equation}
Indeed, if $t\in [1-B, 1]$ and $n\sqrt{1-A}\leq 1$, then
\begin{align*}
n^2(1-|\rho+at+x|)& = n^2(1-A)+ n^2 (A- |\rho+at+x|) \\
& \le
  n^2(1-A)+ n^2 |a| B \leq 1+\f{1+n\sqrt{1-A}} 2\leq 2,
\end{align*}
so that both sides of \eqref{2-2} are bounded up and down by
constant; whereas if $t\in [1-B, 1]$ and $n\sqrt{1-A}\ge 1$, then
\begin{align*}
\left|n\sqrt{1-|\rho+at+x|}-n\sqrt{ 1-A}\right| & \leq \f {n |a||1-t|}{
 \sqrt{1-|at+x+\rho|}+\sqrt{ 1-A}}\\
& \leq \f {n  |a|B} { \sqrt{ 1-A}}
  \leq n^2|a|B= \f  { 1+ n\sqrt{1-A}}2,
\end{align*}
from which \eqref{2-2} follows by triangle inequality. From
($\ref{2-2}$) and ($\ref{2-1}$) with $j=0$, we obtain
\begin{align*}
 \left |\int_{\max\{1-B, -1\}}^1 f(at+x)(1-t)^{\da-1} \xi(t)\,dt \right|
   & \leq c (1+n \sqrt{ 1- A})^{-\mu-\f12}\int_{1-B }^1(1-t)^{\da-1} \, dt\\
& \leq c\, n^{-2\da} |a|^{-\da}(1+n \sqrt{ 1- A})^{-\mu-\f12+\da}.
\end{align*}
If $B \ge \f32$, then the desired inequality \eqref{1-3} follows from the
above inequality. Hence, we assume $B\leq \f32$ from now on.

We now consider the integral over $[-1,1-B]$. Set
$$
\ell= \left \lfloor |\mu|+2\da+\f12 \right \rfloor+1.
$$
Then $1\leq \ell\leq v$ by our assumption.  Since $\xi\in C^\infty [-1,1]$
with $\text{supp}\ \xi\subset [-\f12,1]$, we use ($\ref{2-1}$),
($\ref{2-2}$) and integration by parts $\ell$ times  to obtain
\begin{align*}
&\left | \int_{-1}^{1-B} f(at+x)(1-t)^{\da-1} \xi(t)\, dt\right|
\leq c \sum_{j=1}^\ell |a|^{-j} n^{-2j} (
1+n\sqrt{1-A})^{-\mu-\f12+j} B^{\da-j}\\
& \qquad\qquad\qquad  +c
|a|^{-\ell}\int_{-\f12}^{1-B} |F_\ell
(at+x)|(1-t)^{\da-\ell-1}\, dt\\
&\qquad \le c\, n^{-2\da}| a|^{-\da} (1+n \sqrt{1-A})^{-\mu-\f12+\da}\\
& \qquad \qquad\qquad  +c |a|^{-\ell}n^{-2\ell} \int_{-\f12}^{1-B}(
1+n\sqrt{1-|\rho+x+at|})^{-\mu-\f12+\ell} (1-t)^{\da-\ell-1}\, dt.
\end{align*}
The first term is the desired upper bound in \eqref{1-3}. We only need
to estimate the second term, which we denote by $L$. A change of variable
$s=|a| (1-t)$ shows that
\begin{align*}
 L  &\, := n^{-2\ell} |a|^{-\da} \int_{B |a|}^{\f32 |a|}
    ( 1+ n \sqrt{ 1-|a+x+\rho -s \cdot \text{sgn} \, a|})^{-\mu-\f12+\ell}
   s^{\da-\ell-1}\, ds \\
        & \, =  n^{-2\ell} |a|^{-\da} ( L_1+ L_2)
 \end{align*}
where $L_1$ and $L_2$ are integrals over the intervals
$I_1=\left[ |a|B, \  \f32|a|\right] \cap \left[ 0, \  \f{1-A}2\right]$ and
$I_2=\left[ |a|B, \  \f32|a|\right] \cap \left[ \f{1-A}2,\ \infty\right)$, respectively.
If $s \in I_1$ then
$$
\left|A-   |a+x+\rho -s \cdot \text{sgn}\, a| \right| \le |s| \le
(1-A)/2
$$
so that $ 1-|a+x+\rho -s \cdot \text{sgn} \, a| \sim  1-A$ by triangle
inequality. Consequently,
\begin{align*}
 L_1 & := \int_{I_1} ( 1+ n \sqrt{
 1-|a+x+\rho -s \cdot \text{sgn} \, a|})^{-\mu-\f12+\ell} s^{\da-\ell-1}\, ds\\
& \leq c  ( 1+ n \sqrt{1-A})^{-\mu-\f12+\ell}\int_{B |a|}^{\infty}s^{\da-\ell-1}\,dt\\
& \leq c (1+n \sqrt{ 1-A})^{-\mu-\f12+\ell}(|a|B)^{\da-\ell}\\
& \leq c n^{2\ell-2\da}(1+n \sqrt{ 1- A})^{-\mu-\f12+\delta}.
\end{align*}
If $s \in I_2$, then $s \ge (1-A)/2$ and
$1-|a+x+\rho -s \cdot \text{sgn} \, a|\leq 1-A+s\sim s$ by
triangle inequality. Consequently, since $\ell\ge \mu+\f12$, it follows that
\begin{align*}
L_2&:= \int_{I_2} ( 1+ n
\sqrt{ 1-|a+x+\rho -s \cdot \text{sgn} \, a|})^{-\mu-\f12+\ell} s^{\da-\ell-1}\, ds\\
 & \leq c  \int_{I_2} ( 1+ n \sqrt{s})^{-\mu-\f12+\ell} s^{\da-\ell-1}\,
ds\\
&\leq c n^{-\mu-\f12+\ell}  \int_{ |a|B}^\infty s^{-\f\mu2+\da-\f\ell2-\f54}\, ds,
\end{align*}
since $n^2  |a| B \ge \f 12$. Using the fact that $\ell> -\mu+2\da-\f12$, we
obtain
\begin{align*}
 L_2 \leq c n^{-\mu-\f12+\ell} (|a|B)^{-\f\mu2+\da-\f\ell2-\f14}
 & = c n^{2\ell-2\da}
 (1+n\sqrt{1-A})^{-\f\mu2+\da-\f14-\f\ell2} \\
 &\leq c n^{2\ell-2\da} (1+n\sqrt{1-A})^{-\mu+\da-\f12},
 \end{align*}
using the inequality $\ell\ge \mu+\f12$. Putting these estimates
together completes the proof of \eqref{1-3}.
\end{proof}

\begin{lem} \label{lem:3.5}
Let $\k_j>0$, $a_j\neq 0$, $\xi_j\in C^\infty[-1,1]$ with
$\supp  \xi_j \subset [-\f12,1]$ for $j=1,2,\cdots m$, and
let $\sum_{j=1}^m |a_j| \le 1$. Define
\begin{equation}
 f_m(x):= \int_{[-1,1]^m} P_n^{(\al,\be)} \biggl(\sum_{j=1}^m a_j t_j +x\biggr)
\prod_{j=1}^m \xi_j(t_j) (1-t_j)^{\k_j-1} \, dt
\end{equation}
for $|x| \le1- \sum_{j=1}^m |a_j|$. If $\al\ge \be$, then
\begin{equation}
 |f_m(x)| \leq c  \prod_{j=1}^m  |a_j|^{-\k_j}n^{\al-2\k_j} \left(1+n \sqrt{
1-|A_m+x|}\right)^{-\al-\f12+\tau_m},
\end{equation}
where $A_m:=\sum_{j=1}^m a_j$ and $\tau_m :=\sum_{j=1}^m \k_j$.
\end{lem}

\begin{proof}
Since $ n^{-\al}P_n^{(\al,\be)} (x) \in \mathcal{S}_n^{v_1}(0, 1, \al)$
for $v_1:=\lfloor |\al|+2\k_1\rfloor+4$, we can apply Lemma \ref{lem:3.4}
to conclude that
\begin{align*}
n^{-\a}|f_1(x)| & =
n^{-\al} \left |\int_{-1}^1 P_n^{(\al,\be)}
(a_1t_1+x)(1-t_1)^{\k_1-1} \xi_1(t_1)\, dt_1\right |\\
& \leq c |a_1|^{-\k_1} n^{-2\k_1}
 \left(1+n\sqrt{1-|a_1+x|}\right)^{-\al-\f12+\k_1},
\end{align*}
where $|a_1|+|x|\leq 1$. Hence, the conclusion of the lemma holds
when $m=1$.

Assume that the conclusion of the lemma has been proved for a positive
integer $m$, we now consider the case of $m+1$. Let $v_{m+1}=
\lfloor |\al-\tau_m|+2\k_{m+1}\rfloor +4$. For $i =0,1,\ldots,v_{m+1}$ we define
$$
F_i (x) = C_{n,i}\int_{[-1,1]^m} P_{n+i}^{(\al-i,\be-i)}
 \biggl (\sum_{j=1}^m a_j t_j+x \biggr) \prod_{j=1}^m (1-t_j)^{\k_j-1}\xi_j(t_j)\, dt,
$$
where $C_{n,0} =1$ and $ C_{n,i} =  2^{i}/ \prod_{l=1}^{i} (n+\al+\be+1-l)
 = \CO(n^{-i})$ for $i=1,\ldots, v_{m+1}$.  Using  \eqref{D-Jacobi}, it is easy to
verify that $F_i^{(i)} (x)=f(x)$ for $i=0,1,\cdots, v_{m+1}$. Furthermore, the
induction hypothesis shows that
$$
  |F_i(x)|\leq c \prod_{j=1}^m |a_j|^{-\k_j} n^{\al-2\k_j-2i}
     (1+n\sqrt{1-|A_m+x|})^{-\al-\f12+\tau_m+i}
$$
for $i=0,1,\cdots, v_{m+1}$,  where $|x|+\sum_{j=1}^m|a_j|\leq 1$. By the
definition of $\CS_n^v(\rho,r,\mu)$, this shows that
$$
\prod_{j=1}^m |a_j|^{\k_j} n^{-\al+2\k_j} f_m(x) \, \in \,
  \mathcal{S}_n^{v_{m+1}} \biggl( A_m,  1-\sum_{j=1}^m|a_j|, \al-\tau_m \biggr).
$$
Since $v_{m+1} \ge |\al-\tau_m|+2\k_{m+1} +\f32$, we can then apply Lemma
\ref{lem:3.4} to the integral
$$
 f_{m+1}(x) = \int_{-1}^1 f_m(a_{m+1}
     t_{m+1} +x)(1-t_{m+1})^{\k_{m+1}-1}\xi_{m+1}(t_{m+1})\, dt_{m+1}
$$
to conclude that
\begin{align*}
 &  \prod_{j=1}^m |a_j|^{\k_j} n^{-\al+2\k_j} \left| f_{m+1}(x)\right| \\
 &  \qquad\qquad   \leq c n^{-2\k_{m+1}} |a_{m+1}|^{-\k_{m+1}}
       (1+n\sqrt{1-|A_{m+1}+x|})^{-\al-\f12 +\tau_{m+1}},
\end{align*}
where $ |x|+|a_{m+1}|\leq 1 -\sum_{j=1}^m |a_j|$. This completes the induction
and the proof.
\end{proof}

We are now in a position to prove Theorem \ref{MainEst}

\medskip\noindent
{\it Proof of Theorem \ref{MainEst}}.  Let $\psi\in C^\infty[-1,1]$
satisfy $\psi(t)=1$ for $\f12\leq t\leq 1$, and $\psi(t)=0$ for
$-1\leq t \leq -\f12$.  We define
\begin{align*}
\begin{split}
\xi_{1,j}(t)&=\vi_j(t) \psi(t) (1+t)^{\k_j-1},\\
\xi_{-1, j} (t) &= \vi_j(-t) (1-\psi(-t))(1+t)^{\k_j-1},
\end{split} \qquad\quad
j=1,\cdots, m.
\end{align*}
Evidently, $\xi_{1, j}, \xi_{-1,j}\in C^\infty [-1,1]$ and
$\text{supp}\, \xi_{1, j}, \text{supp}\, \xi_{-1, j}\subset [-\f12, 1]$. Since
\begin{align*}
&\int_{-1}^1 g(t_j)\vi_j(t_j)(1-t_j^2)^{\k_j-1}\, dt_j \\
& \qquad  = \int_{-1}^1 g(t_j)\xi_{1,j}(t_j)(1-t_j)^{\k_j-1}\, dt_j
     + \int_{-1}^1 g(-t_j)\xi_{-1,j}(t_j)(1-t_j)^{\k_j-1}\, dt_j,
\end{align*}
we can write
\begin{align*}
J: &=  \int_{[-1,1]^m} P_n^{(\al,\be)}\biggl(\sum_{j=1}^m
a_jt_j+x\biggr)\prod_{j=1}^m \vi_j(t_j)(1-t_j^2)^{\k_j-1}\, dt \\
& = \sum_{\ve \in \{1,-1\}^m} \int_{[-1,1]^m}
P_n^{(\al,\be)}\biggl(\sum_{j=1}^m \ve_ja_jt_j+x\biggr)
\prod_{j=1}^m \xi_{\ve_j,j}( t_j)(1-t_j)^{\k_j-1}\, dt\\
 & =: \sum_{\ve \in \{1,-1\}^m} I_{\ve}(x).
\end{align*}
Recall $|a| =  \sum_{j=1}^m |a_j|$. For $\ve \in \{1,-1\}^m$, we write
$a(\ve) : = \sum_{j=1}^m a_j\va_j$. Applying Lemma \ref{lem:3.5}
to $I_\ve$ gives
\begin{align*}
  |I_\ve(x)| & \leq c n^{\a -2|\k|} \prod_{j=1}^m |a_j|^{-\k_j}
     \left(1+ n \sqrt{ 1-|x+ a(\ve)|}\right)^{-\al-\f12 +|\k|}\\
   & \leq n^{\al-2|\k|} \prod_{j=1}^m|a_j|^{-\k_j} \left(1+ n \sqrt{
        1-|x|- |a|}\right)^{-\al-\f12 +|\k|}
\end{align*}
for each $\ve \in \{1,-1\}^m$, where we have used the assumption
$\al\ge |\k|-\f12$ and the inequality $|x+ \sum_{j=1}^m \va_j a_j|
\leq |x|+\sum_{j=1}^m |a_j|$ in the last step. Consequently,
\begin{align}\label{3-9-07}
 |J| \leq c\, 2^m  n^{\al-2|\k|} \prod_{j=1}^m |a_j|^{-\k_j} \left(1+ n
\sqrt{ 1-|x|- |a|}\right)^{-\al-\f12 +|\k|}.
\end{align}

Finally, we claim that the desired inequality \eqref{eq:MainEst}
is a consequence of (\ref{3-9-07}). In fact,  without loss of
generality, we may assume that
\begin{equation}\label{3-10-07}|a_j|\ge
  n^{-1}\sqrt{1-|a|-|x|}+n^{-2},\quad \text{for $j=1,\cdots,p$}
 \end{equation}
and
\begin{equation}\label{3-11-07} |a_j|< n^{-1}\sqrt{1-|a|-|x|}+n^{-2},\quad
 \text{for $j=p+1,\cdots,m$}.
\end{equation}
We then apply (\ref{3-9-07}) with $m$ and $x$
replaced by $p$ and $\sum_{j=p+1}^m a_j t_j +x$, respectively, to
obtain
\begin{align*}
M_p(x,t'):= & \left |\int_{[-1, 1]^p} P_n^{(\a,\b)}\biggl(\sum_{j=1}^m a_j t_j+x\biggr)
      \prod_{j=1}^p \vi_j(t_j) (1-t_j^2)^{\k_j-1}\, dt \right| \\
 \leq & c\, 2^p n^{\al-2\sum_{j=1}^p\k_j} \prod_{j=1}^p
|a_j|^{-\k_j} \left(1+n A(x) \right)^{-\al-\f12+\sum_{j=1}^p
\k_j},
\end{align*}
where $t':= (t_{p+1}, \cdots, t_m) \in [-1,1]^{m-p}$ and
$A(x):=\sqrt{1-|a|-|x|}$, and we have used the inequality
$|\sum_{j=p+1}^m a_j t_j +x|\leq \sum_{j={p+1}}^m|a_j|+|x|$ as well
as the fact that $\al \ge \sum_{j=1}^p \k_j -\f12$.
Using the assumption ($\ref{3-10-07}$), we then
obtain
\begin{align*}
M_p(x,t')& \leq  c\, n^{\al-2|\k|} \prod_{i=1}^p |a_i|^{-\k_i}
\prod_{j={p+1}}^m
(n^{-1}A(x)+n^{-2})^{-\k_j}  \left(1+n A(x)\right)^{-\al-\f12+|\k|}\\
& \le c n^{\al-2|\k|} \prod_{j=1}^m
\left(|a_j|+n^{-1}A(x)+n^{-2}\right)^{-\k_j}  \left(1+n A(x) \right)^{-\al-\f12+|\k|}.
\end{align*}
Consequently, it follows that
\begin{align*}
& \left |\int_{[-1, 1]^m} P_n^{(\a,\b)}\biggl(\sum_{j=1}^m a_j t_j+x \biggr)
   \prod_{j=1}^m \vi_j(t_j) (1-t_j^2)^{\k_j-1}\, dt \right | \\
& \leq \int_{[-1,1]^{m-p}} M_p(x,t')
  \prod_{i=p+1}^m \vi_i(t_i) (1-t_i^2)^{\k_i-1}\, dt_{p+1}\cdots dt_m\\
&\leq c  n^{\al-2|\k|} \prod_{j=1}^m
\left(|a_j|+n^{-1}A(x)+n^{-2}\right)^{-\k_j}  \left (1+n A(x)\right)^{-\al-\f12+|\k|},
\end{align*}
proving the desired inequality \eqref{eq:MainEst}. \qed


\subsection{Proof of the pointwise estimate of the kernel on the sphere}

For estimating the kernel, we will need information on the $(C,\d)$ means
of the Jacobi expansion. We start with a  result in \cite[p. 261, (9.4.13)]{Szego}
and its extension in \cite{Li} given in the following lemma.

\begin{lem} \label{lem:3.6}
For any $\alpha,\beta > -1$ such that $\alpha+ \beta+ \delta + 3 >0$,
$$
  K_n^\delta(w^{(\alpha,\beta)}, 1,u) = \sum_{j=0}^J b_j(\alpha,\beta,\delta,n)
    P_n^{(\alpha+\delta+j+1,\beta)}(u) + G_n^{\delta}(u),
$$
where $J$ is a fixed integer and
$$
   G_n^\delta(u) = \sum_{j=J+1}^\infty d_j(\alpha,\beta,\delta,n)
  K_n^{\delta+j}(w^{(\alpha,\beta)}, 1,u);
$$
moreover, the coefficients satisfy the inequalities,
$$
 |b_j(\alpha,\beta,\delta,n)| \le c n^{\alpha+1-\delta - j} \quad \hbox{and}
 \quad  |d_j(\alpha,\beta,\delta,n)| \le c j^{-\alpha-\beta-\delta - 4}.
$$
\end{lem}

Since the kernel function $K_n^{\delta + j}(w^{(\alpha,\beta)}, 1,u)$
contained in the $G_n^\delta$ term has larger index, it could be handled
by using the following estimate of the kernel function, which was
used in \cite{BC} and \cite{CTW} (see Theorem 3.9 there).

\begin{lem} \label{lem:3.7}
Let $\alpha,\beta \ge -1/2$.  If $\delta \ge \alpha + \beta + 2$, then
$$
 |K_n^\delta(w^{(\alpha,\beta)}, 1,u)| \le c n^{-1}
        (1-u+n^{-2})^{-( \alpha+3/2)}.
$$
\end{lem}

\medskip
\noindent
{\it Proof of Theorem \ref{thm:S-estimate}.}
We start from the integral expression \eqref{kernel} of $K_n^\d(h_\k^2;x,y)$.
The first step of the proof is to replace the kernel
$K_n^\d(w^{(\l_\k-\frac12,\l_\k-\frac12)})$ by
the expansion in Lemma \ref{lem:3.6}.  Let $\a= \b = |\k| + (d-2)/2$ and let
$J = \lfloor \a+\b+2 \rfloor= \lfloor 2 |\k| + d \rfloor$. The choice of $J$
guarantees that we can apply Lemma \ref{lem:3.7} on $G_n^\d$ term.
Combining the formula \eqref{kernel} and Lemma \ref{lem:3.6}, we obtain
\begin{equation*} 
 K_n^\delta(h_\kappa^2; x,y) = \sum_{j=0}^J b_j(\alpha,\beta, \delta,n)
    \Omega_j(x,y) + \Omega_*(x,y),
\end{equation*}
where
\begin{equation*} \label{eq:3.5}
\Omega_j(x,y) =  c_\kappa \int_{[-1,1]^{d+1}} P_n^{(\a+\d + j + 1,\b)}(u(x, y, t))
 \prod_{i=1}^{d+1} (1+t_i)(1-t_i^2)^{\kappa_i-1} d t,
\end{equation*}
and
\begin{equation*} \label{eq:3.6}
\Omega_*(x,y) =  c_\kappa \int_{[-1,1]^{d+1}} G_n^\delta(u(x,y,t))
 \prod_{i=1}^{d+1} (1+t_i)(1-t_i^2)^{\kappa_i-1} d t,
\end{equation*}
in which $u(x,y,t) = x_1y_1t_1 + \ldots + x_{d+1} y_{d+1} t_{d+1}$.

Since the index of the Jacobi polynomial in $\Omega_0$ are $\a + \d + 1 =
\d + |\k| + \frac{d}{2}$ and $|\k| +\frac{d-2}{2}$, we can use
Theorem \ref{MainEst} with $m = d+1$, $x =0$ and $a_j = x_j y_j$ to
estimate $\Omega_0$ for all $\d > -1$. Using the fact that
$1 - \la \bar x, \bar y\ra = \|\bar x -\bar y\| /2$ for $x ,y \in S^d$, this
shows that $b_0(\a,\b,\d,n)\Omega_0$ is bounded by the first term in
the right hand of \eqref{eq:MainEst}. The same estimate evidently holds
for $\Omega_j$. The estimate of $\Omega_*$ uses Lemma \ref{lem:3.7},
which can be handled easily as shown in \cite{LX}.

Finally, we note that Theorem \ref{MainEst} can also be applied to
the kernel $P_n(h_\k; x,y)$ in \eqref{proj-Kernel}, which gives the
pointwise estimate of \eqref{eq:S-estP}.
\qed


\subsection{Proof of the pointwise estimate of the kernel on the simplex}

Recall the formula for $\Kb_n^\d(W_\k^T; x,y)$ in \eqref{kernelT}.
Setting $\a = |\k|+(d-2)/2$ and $J = \lfloor \a+3/2 \rfloor$, we again use
Lemma \ref{lem:3.6} to break the kernel $\Kb_n^\d (W_\k^T;x,y)$ into a sum
$$
\Kb_n^{\delta}(W_{\k}^T; x,y)
    =\sum_{j=0}^J b_j(\a,-1/2,\d,n)\Omega_j(x,y) +\Omega_{\ast}(x,y),
$$
where
\begin{align*}
\Omega_j(x,y)=c_{\k}
  \int_{[-1,1]^{d+1}}P_n^{(\a +\d+j+1,- \frac{1}{2})}(2 z(x,y,t)^2-1)
 \prod_{i=1}^{d+1}(1-t_i^2)^{\k_i-1} dt
\end{align*}
and
$$
\Omega_*(x,y) =  c_\kappa \int_{[-1,1]^{d+1}} G_n^\delta(2z(x, y,
t) ^2-1)
  \prod_{i=1}^{d+1}(1-t_i^2)^{\k_i-1}d t,
$$
where $z(x, y, t)=\sum_{j=1}^{d+1} \sqrt{x_jy_j} t_j $,
$x_{d+1}=\sqrt{1-|x|}$ and $y_{d+1} =\sqrt{1-|y|}$. Using the
quadratic transform $P_n^{(\l, -\f12)} (2t^2-1)=a_n
P_{2n}^{(\l,\l)}(t)$ with $a_n =O(1)$, we have
$$\Omega_j(x,y)=O(1)\int_{[-1,1]^{d+1}}P_{2n}^{(\a +\d+j+1, \a +\d+j+1)}(z(x, y,
t)) \prod_{i=1}^{d+1}(1-t_i^2)^{\k_i-1} dt.$$
 Since $\xi,\zeta \in S^d$, we have $1 - z(x,y,t) \ge \|\xi
- \zeta\|^2 /2$. Hence, we can follow the same procedure as in the
proof of Theorem \ref{thm:S-estimate} to prove Theorem
\ref{thm:T-estimate}.


\section{Lower bound estimate}
\setcounter{equation}{0}

The lower bound estimate comes down to the proof of Proposition
\ref{prop:Tn-lwbd}, which gives a lower bound of
$T_n^\d(w_{\l,\mu})$ in \eqref{eq:Tn-d} for $\d \le \mu$. The case
of $\d = \mu$ is already established in \cite{LX}, but the proof
there is rather involved and may not work for the case $\d < \mu$.
Below we shall follow a different and simpler approach, which
works for $\d \le \mu$ and gives, in particular, a simpler proof
in the case of $\d = \mu$.

\medskip\noindent
{\it Proof of Proposition \ref{prop:Tn-lwbd}.}
It is known that
\begin{align*}
 & T_n^\d(w_{\l,\mu}; 1) \ge  T_n^\d(w_{\mu,\l}; 0) = c n^{\l+\mu-\d+ \frac{1}{2}} \\
 & \quad \times \int_{0}^1 \Big|
\int_{-1}^1 P_n^{(\lambda+ \mu +\d+\frac12, \lambda+\mu -\frac12)} (s t)
(1-s^2)^{\mu-1} ds\Big| t^{2\mu}(1-t^2)^{\l -\frac12} dt + \CO(1).
\end{align*}
This is proved in \cite[p. 293]{LX}, where the equation is stated
for $T_n^\d(w_{\l,\mu};0)$ and we should mention that in the last
two displayed equations in \cite[p. 293]{LX}, $w_{\mu,\l}$ should
have been $w_{\l,\mu}$. As  a result of this relation, we see that
Proposition \ref{prop:Tn-lwbd} follows from the lower bound of the
double integral of the Jacobi polynomial given in the next
proposition. \qed

\medskip

\begin{prop} \label{prop:lower}
Assume $ \l, \mu\ge 0$ and $ \l \ge \d >-1$. Let $a = \l + \mu +
\d $ and $b = \l+\mu -1$. Then
\begin{align} \label{lowerJacobi}
&\int_{0}^1 \left| \int_{-1}^1 P_n^{(a+\frac12,b+\frac12)} (t y)
(1-t^2)^{\mu-1} dt \right| |y|^{2\mu} (1-y^2)^{\l -1/2} dy \\
 & \qquad\qquad \ge c \; n^{-\mu -1/2}  \begin{cases}  \log n, & \text{if $\d = \l$},\\
 1, &\text{if $ -1<\d < \l$},
  \end {cases}  \notag
\end{align}
where, when $\mu =0$, the inner integral is defined in the sense of
(\ref{1-8-07}).
\end{prop}

Let us denote the left hand side of
\eqref{lowerJacobi} by $I_n$. First, we assume that $0<\mu< 1$.
Changing variables $t = u /y$, followed by $y = \cos \phi$ and $u
= \cos \t$, and restricting the range of the outside integral lead
to
\begin{align*}
I_n & \ge c \int_{\sqrt{2}/2}^{1} \left |
   \int_{-y}^y P_n^{(a+\frac12,b+\frac12)} (u) y
 (y^2-u^2)^{\mu-1} du \right | (1-y^2)^{\l -1/2} dy \\
& \ge c \int_{ n^{-1}}^{\pi/4}
   \left | \int_{\phi}^{\pi -\phi} P_n^{(a+\frac12,b+\frac12)} (\cos \t)
           (\cos^2 \phi-\cos^2 \t)^{\mu-1} \sin \t d \t \right | (\sin \phi)^{2\l}  d\phi.
\end{align*}
We need the asymptotics of the Jacobi polynomials as given in
\cite[p. 198]{Szego},
$$
 P_n^{(\alpha,\beta)}(\cos \theta) =
  \pi^{-\frac12} n^{-\frac12}   (\sin \tfrac{\t}2 )^{-\a-\frac12}
  (\cos \tfrac{\t}2)^{-\b-\frac12}
        \left[ \cos (N\theta+\tau) + \CO(1) (n \sin\theta)^{-1} \right]
$$
for $ n^{-1} \le \t \le \pi -  n^{-1}$, where $N=
n+\frac{\alpha+\beta+1}2$ and $\tau = -\frac{\pi}2
({\alpha+\frac12}).$ Applying this asymptotic formula with $\a =
a+1/2$ and $\b= b+1/2$ we obtain
\begin{equation} \label{eq:lower1}
  I _n \ge  c \,n^{-1/2} \int_{ n^{-1}}^{\pi/4} |M_n(\phi)| (\sin \phi)^{2\l}
       d\phi  -    \CO(1)  E_n,
\end{equation}
where $M_n(\phi)$ is the integral over the main term of the asymptotics
\begin{equation}\label{4-3-07}
   M_n(\phi): = \int_\phi^{\pi-\phi}  \frac{(\cos^2 \phi-\cos^2 \t)^{\mu-1} }
        { (\sin \frac{\t}2 )^{a } (\cos \frac{\t}2)^{b} } \cos (N\t + \tau)
        d\t,
\end{equation}
and $E_n$ comes from the remainder term in the asymptotics
\begin{equation}\label{4-4-07}
 E_n : = n^{-\frac32}   \int_{ n^{-1}}^{\pi/4}
    \int_\phi^{\pi-\phi}  \frac{(\cos^2 \phi-\cos^2 \t)^{\mu-1} }
        { (\sin \frac{\t}2 )^{a+1} (\cos \frac{\t}2)^{b+1} } d\t
          \, (\sin \phi)^{2\l}  d\phi.
\end{equation} Here $ N=n+\f{a+b}2+1$ and $\tau= -\f\pi2(a+1)$.

In order to handle the main part of \eqref{eq:lower1}, we first
derive an asymptotic formula for $M_n(\phi)$. We need the
following lemma, which follows directly from \cite[p.
49]{A.Erdelyi}.

\begin{lem} \label{Erdelyi}
If $0<\mu< 1$,  $g(t)$ is continuously differentiable on the
interval $[\al,\b]$,  and $\xi\in\mathbb{R}-\{0\}$ then
\begin{align*}
  &\int_\a^\b g(t)  e^{i \xi  t} (t-\a)^{\mu-1} (\b-t)^{\mu-1} dt \\
  & \qquad  = \Gamma(\mu) |\xi|^{-\mu}  \left[ e^{- \f{i \pi \mu\xi}{2|\xi|}} g(\b) e^{i \xi \b} +
                 e^{\f{i \pi \mu\xi}{2|\xi|}}g(\a)  e^{i \xi \a} \right]+ R_\xi,
\end{align*}
as $|\xi| \to + \infty$, where
$$
    |R_\xi| \le  |\xi|^{-1} \int_\a^\b |g'(t) | (t-\a)^{\mu-1} (\b -t)^{\mu-1} dt.
$$
\end{lem}

\begin{lem} \label{LowerLem1}
Assume  $0<\mu< 1$, $\l\ge 0$ and $\l\ge \da\ge  -1$. Let
$M_n(\phi)$ be defined by \eqref{4-3-07}. Then for $ 0< \phi \le
\pi/4$,
\begin{equation}\label{4-5-07}
M_n(\phi) =K_n(\phi) +G_n(\phi),\end{equation} where
\begin{align}
K_n(\phi) = \,  & \Gamma(\mu)  N^{- \mu}
 \frac{2^a(\sin (2\phi))^{\mu-1} }{(\pi - 2\phi)^{\mu-1} (\sin \phi)^a }\notag \\
   &   \times \left[ (-1)^n (\sin \tfrac{\phi}2)^{a-b}
   \cos\bigl( N\phi +\g +\tfrac{(a-b)\pi}2\bigr) +
      (\cos \tfrac{\phi}2)^{a-b} \cos (N\phi +\g) \right],\label{4-6-07}\end{align}
$\g=\tau+\f{\pi\mu}2$, and the remainder satisfies
\begin{equation}\label{Rn-bound}
    |G_n(\phi)| \le c n^{-1}
    \phi^{\mu-\l-\da-2}.\end{equation}
\end{lem}

\begin{proof}
Writing $\cos (N \t + \tau) =( e^{i (N \t + \tau)} + e^{-i (N \t + \tau)} )/2$, we
split $M_n(\phi)$ into two parts, $M_n^+(\phi)$ and $M_n^- (\phi)$, respectively,
and apply Lemma \ref{Erdelyi} to these integrals. For $M_n^+(\phi)$
we define a function $f_\phi$ as
$$
  f_\phi(\theta) =  \frac{(\cos^2 \phi-\cos^2 \t)^{\mu-1} }
        { (\sin \frac{\t}2)^a (\cos \frac{\t}2)^{b}  (\t-\phi )^{\mu-1}
            (\pi -\phi - \t)^{\mu-1}}
$$
for $\phi <  \t < \pi -\phi$ and define its value at the boundary
by limit. Then it is easily seen  that
\begin{equation*}
f_\phi(\ta) =\begin{cases} \left(\df{
\sin(\pi-\phi-\ta)\sin(\ta-\phi)}{(\pi-\phi-\ta)(\ta-\phi)}\right)^{\mu-1}
\df 1{(\sin \f\ta 2)^a (\cos\f\ta2)^b },\   \  &\text{if
$\ta\in(\phi,\pi-\phi)$},\\[5mm]
\left( \df{\sin(\pi-2\phi)}{\pi-2\phi}\right)^{\mu-1} \df1{ (\sin
\f\ta2)^a (\cos\f\ta2)^b},\   \ &\text{if $\ta=\phi$ or
$\pi-\phi$, }\end{cases}
\end{equation*}
 is continuously
differentiable on $[\phi,\pi-\phi]$. Hence, invoking  Lemma
\ref{Erdelyi} with $\xi=N$,  and by  a straightforward
computation, we obtain
\begin{align*}
 M_n^+(\phi)  & =  \frac{e^{i\tau}}{2}
    \int_\phi^{\pi-\phi} f_\phi(\t) e^{i N\t}(\t-\phi)^{\mu-1}  (\pi -\phi - \t)^{\mu-1} d\t \\
& = \Gamma(\mu)  N^{-\mu}
 \frac{(\sin (2\phi))^{\mu-1} }{(\pi - 2\phi)^{\mu-1}}
 \f{2^{a-1}} {(\sin{\phi})^a } \\
& \qquad
   \times \displaystyle\left[ (\sin \tfrac{\phi}2)^{a-b}
   e^{i [N(\pi-\phi)-\f{\pi\mu}2+\tau]} + (\cos \tfrac{\phi}2)^{a-b}
      e^{i[N\phi+\f{\pi\mu}2+\tau]} \right]  +  R_n^+(\phi),
\end{align*}
in which
$$
|R_n^+(\phi)| \le N^{-1} \int_{\phi}^{\pi -\phi} |f'_\phi(\t)| (\t -\phi)^{\mu-1}
     ( \pi - \phi -\t)^{\mu-1} d\t.
$$
Since $0<\phi\leq \f\pi4$, using the fact that $\sin x /x$ is analytic and
that $\sin (\pi - \t -\phi) = \sin (\t + \phi)$, from the definition
of $f_\phi$ we see easily that for $\phi<\ta<\pi-\phi$,
$$|f_\phi'(\ta)|\leq c (\ta^{\mu-a-2}+(\pi-\ta)^{\mu-b-2}).$$
 This implies that for  $0<\phi\leq \f\pi4$,
\begin{align*}
 |R_n^+(\phi)|
  & \le c N^{-1} \left[ \int_{\phi}^{\pi/2} \t^{\mu-a-2} (\t -\phi)^{\mu-1} d \t
  +  \int^{\pi-\phi}_{\pi/2}  (\pi-\t)^{\mu-b-2}(\pi-\t-\phi)^{\mu-1} d\t \right ]\notag\\
  &\leq c n^{-1}  \int_{\phi}^{\pi/2} \t^{\mu-a-2} (\t -\phi)^{\mu-1} d \t\notag
\end{align*}
as $a \ge b$ and the the first term dominates.  A simple
computation shows then
\begin{align}\label{Rn-bound-2}
 |R_n^+(\phi)|   &\leq c
  n^{-1} \phi^{\mu-a-2} \int_{\phi}^{2\phi}  (\t -\phi)^{\mu-1} d \t
  + c n^{-1}  \int_{2\phi}^{\pi/2} \ta^{2\mu-a-3}  d \t \\
  &\leq
  c n^{-1} \phi^{2\mu-2-a}= cn^{-1} \phi^{\mu-\l-\da-2}\notag\end{align}
  since $a=\l+\mu+\da > 2\mu-2$.

Similarly,  using  Lemma \ref{Erdelyi} with $\xi=-N$, we  derive a
similar relation for $M_n^-(\phi)$:
\begin{align*}
 M_n^-(\phi)  & =  \frac{e^{-i\tau}}{2}
    \int_\phi^{\pi-\phi} f_\phi(\t)
     e^{-i N\t}(\t-\phi)^{\mu-1}  (\pi -\phi - \t)^{\mu-1} d\t \\
& = \Gamma(\mu)  N^{-\mu}
 \frac{(\sin (2\phi))^{\mu-1} }{(\pi - 2\phi)^{\mu-1}}
 \f{2^{a-1}} {(\sin{\phi})^a } \\
& \qquad
   \times \left[ (\sin \tfrac{\phi}2)^{a-b}
   e^{-i [N(\pi-\phi)-\f{\pi\mu}2+\tau]} + (\cos \tfrac{\phi}2)^{a-b}
      e^{-i[N\phi+\f{\pi\mu}2+\tau]} \right]  +  R_n^-(\phi),
\end{align*}
where the  error term $R_n^-(\phi)$ satisfies the same upper bound
as in \eqref{Rn-bound-2}. Since $M_n (\phi) = M_n^+(\phi) +
M_n^-(\phi)$ and $N\pi+2\tau=n\pi+\f{b-a}2\pi$,  the desired
expression for $M_n(\phi)$ follows with $G_n(\phi) = R_n^+(\phi) +
R_n^-(\phi)$, which satisfies the stated bound.
\end{proof}

\begin{lem}\label{4-4-lem}
Assume that  $0<\mu< 1$, $\l\ge 0$ and $\l \ge \da >  -1$. Then
$$
\int_{n^{-1}}^{\f \pi 4} |M_n (\phi)|(\sin\phi)^{2\l}\, d\phi
\ge  c n^{-\mu}
\begin{cases} \log n,&\  \  \text{if $\l=\da$},\\
 1,&\  \  \text{if $-1<\da<\l$}.
\end{cases}
$$
\end{lem}

\begin{proof}
Since  $a-b=\da+1>0$,  we can  choose an absolute constant $\va\in
(0, \f \pi4)$ satisfying  $(\tan \f \va 2)^{a-b} \leq \f14$. We
then use  ($\ref{4-6-07}$),  and obtain that  for $\phi\in (0,
\va)$,
\begin{align*}
|K_n(\phi)|& \ge c n^{-\mu} \phi^{-\l-\da-1} \left( |\cos
(N\phi+\g)|-
    \left(\tan \tfrac{\phi}2\right)^{a-b}\right)\\
&\ge c n^{-\mu} \phi^{-\l-\da-1} \left( \cos^2 (N\phi +\g)- \f14\right)\\
&= \f c4n^{-\mu} \phi^{-\l -\da-1} + \f c2 n^{-\mu} \phi^{-\l
-\da-1}\cos
  (2N\phi +2\g),
\end{align*}
where we have used the fact that $(\tan \f\phi2)^{a-b} \leq (\tan
\f\va 2)^{a-b}\leq  \f14$  for $0<\phi\leq \va$ in the second
step,  and the identity $\cos^2t =\f12 +\f12\cos 2t$ in the last
step. It follows that
\begin{align*}
\int_{n^{-1}}^{\va} |K_n (\phi)|(\sin\phi)^{2\l}\, d\phi&\ge
cn^{-\mu} \int_{n^{-1}}^{\va}\phi^{\l-\da-1}\, d\phi \\
&\qquad  +c
n^{-\mu}\int_{n^{-1}}^{\va}\phi^{\l-\da -1} \cos (2N\phi +2\g)\, d\phi \\
&\ge  c n^{-\mu}
\begin{cases} \log n,&\  \  \text{if $\l=\da$},\\
 1,&\  \  \text{if $-1<\da<\l$}.
\end{cases}
\end{align*}
where we have used an integration by parts in the last step.

To complete the proof, we just need to  observe that by
($\ref{4-5-07}$),
$$\int_{n^{-1}}^{\f\pi4} |M_n(\phi)|(\sin\phi)^{2\l} \, d\phi \ge \int_{n^{-1}}^{\va} |K_n
(\phi)|(\sin\phi)^{2\l}\, d\phi-\int_{n^{-1}}^{\va} |G_n
(\phi)|(\sin\phi)^{2\l}\, d\phi,$$ whereas  by ($\ref{Rn-bound}$),
$$
\int_{n^{-1}}^{\va} |G_n (\phi)|(\sin\phi)^{2\l}\, d\phi\leq c
n^{-1}\log n + c n^{-\mu+\da-\l}
$$
which is small than the bound for the first term in magnitude as $0< \mu < 1$.
\end{proof}

\begin{lem}\label{4-5-lem}
 Assume  $0<\mu< 1$, $\l \ge 0$ and $\l \ge \da > -1$.
 Let $E_n$ be defined by \eqref{4-4-07}. Then
 $$
      E_n \leq  c n^{-\mu-\f12-(\l -\da)}+cn^{-\f32} \log n.
 $$
\end{lem}

\begin{proof}
 By  ($\ref{4-4-07}$) and the identity $\cos^2\ta
-\cos^2 \phi=\sin (\ta+\phi)\sin (\t-\phi)$, we obtain
\begin{align*}
E_n&=n^{-\f32} \int_{n^{-1}}^{\f\pi4} \int_{\phi}^{\pi-\phi} \df {
\sin^{\mu-1} (\ta+\phi)\sin^{\mu-1}(\ta-\phi)}{(\sin^{a+1}\f
\t2)(\cos^{b+1}\f\t2)}\, d\t \sin^{2\l}\phi\, d\phi\\
&\leq cn^{-\f32} \int_{n^{-1}}^{\f\pi4} \int_{\phi}^{\f\pi2}
 \t^{\mu-a-2} (\ta-\phi)^{\mu-1}\, d\t \phi^{2\l}\,
 d\phi.
\end{align*}
The inner integral can be estimated by splitting the integral as two
parts, over $[\phi, 2 \phi]$ and over $[2\phi, \pi/2]$, respectively.
Upon considering the various cases and taking into the account that
$a = \l  + \mu + \d$ and $\l \ge 0$, we conclude that
 \begin{align*}
 E_n &\leq c n^{-\f32} \int_{n^{-1}}^{\f\pi4}  \Bl(\phi^{2\l} |\log \phi|
 + \phi^{ \l+\mu-\da-2}\Br)\, d\phi\leq c
 n^{-\mu-\f12-(\l-\da)}+cn^{-\f32} \log n.
\end{align*}
This completes the proof of Lemma \ref{4-5-lem}.
\end{proof}

We now  return to the  proof of Proposition \ref{prop:lower}.

\medskip\noindent
{\it Proof of Proposition \ref{prop:lower} (Continue).}\  \ We
consider the following cases:

\medskip\noindent
{\it Case 1. $0 < \mu <1$.}
This  case  follows directly from
($\ref{eq:lower1}$) and Lemmas $\ref{4-4-lem}$ and $\ref{4-5-lem}$.

\medskip\noindent
{\it Case 2. $\mu=0$ or $1$.} In the case $\mu=0$,  $I_n$ in limit form reduces to
\begin{align*}
    I_n  & = \int_0^1 \left | P_n^{(a+ \frac12, b+\frac12)}(y)
   + P_n^{(a+\frac12, b+\frac12)}(-y) \right | (1-y^2)^{\l-1/2} dy \\
   & \ge  \int_{n^{-1}}^{\pi/4} \left | P_n^{(a+ \frac12, b+\frac12)}(\cos \phi)
   + P_n^{(a+\frac12, b+\frac12)}(\cos (\pi-\phi)) \right | (\sin \phi)^{2\l} d\phi.
\end{align*}
The asymptotic formula of the Jacobi polynomial gives
\begin{align*}
& P_n^{(a+ \frac12, b+\frac12)}(\cos \phi) + P_n^{(a+\frac12, b+\frac12)}(\cos(\pi- \phi))
   = \frac{\pi ^{-1/2} n^{-1/2} }{(\sin \frac{\phi}{2})^{a+1}   (\cos \frac{\phi}{2})^{a+1}} \\
& \quad \times \left[ (\cos \tfrac{\phi}{2})^{a-b} \cos (N\phi+\tau) +
      (\sin \tfrac{\phi}{2})^{a-b}   \cos (N (\pi-\phi)+\tau) \right] +
          \CO\left((n \sin \phi)^{-1}\right),
\end{align*}
which is essentially the same as the asymptotic formula for
$M_n(\phi)$ in Lemma \ref{LowerLem1} with $\mu =0$ and a smaller
remainder.  Thus, a  proof almost identical to that of Lemma
\ref{4-4-lem}  will yield  Proposition \ref{prop:lower} for $\mu=0$.  Proposition \ref{prop:lower} for $\mu=1$ can be proved in a similar way.

\medskip\noindent
{\it Case 3. } $\mu > 1$.  In this case, we  denote by $r$ the largest integer
smaller than $\mu$.  We  then use
\eqref{D-Jacobi} and integrate  by parts $r$ times to obtain
\begin{align*}
\int_{-y}^y & P_n^{(a+\frac12,b+\frac12)} (u) (y^2-u^2)^{\mu -1} du \\
 =& \frac{(-2)^r} {\prod_{i=1}^r(n+a+b+2-i)} \int_{-y}^y
 P_{n+r}^{(a+\frac12-r,b+\frac12-r)} (u)
  \frac{d^r}{du^r} \Big[(y^2-u^2)^{\mu-1} \Big]du.
\end{align*}
Since  $ [(y^2-u^2)^{\mu-1}]^{(r)} = A q(y, u)  (y^2 -u^2)^{\mu - r -1}$,
where $A$ is a  nonzero constant and
 $q(y, u)$ is a polynomial in $y$ and $u$ which satisfies
 $q(y, y)=(-1)^r q(y, -y) =1$, we conclude that
\begin{align*}
& \left | \int_{-y}^y  P_n^{(a+\frac12,b+\frac12)}(u)
(y^2-u^2)^{\mu-1} du \right | \\
& \qquad     \ge \, c  n^{-r} \left | \int_{-y}^y
 P_{n+r}^{(a'+\frac12, b'+\frac12)} (u)
   q(y, u) (y^2-u^2)^{\mu' -1} du   \right |,
 \end{align*}
 where $\mu'=\mu-r\in (0, 1]$, $a'=\l+\mu'+\da$ and
 $b'=\l+\mu'+\da$. It follows that
 \begin{align*}
I_n & \ge c n^{-r}  \int_{\sqrt{2}/2}^{1} \left |
   \int_{-y}^y P_{n+r}^{(a'+\frac12,b'+\frac12)} (u) q(y, u) y
 (y^2-u^2)^{\mu'-1} du \right | (1-y^2)^{\l -1/2} dy \\
& \ge c n^{-r} \int_{ n^{-1}}^{\pi/4}
   \left | \int_{\phi}^{\pi -\phi} P_{n+r}^{(a'+\frac12,b'+\frac12)} (\cos
   \t)q_\phi (\cos\ta)
           (\cos^2 \phi-\cos^2 \t)^{\mu'-1} \sin \t d \t \right | \\
            & \qquad\qquad\qquad \times(\sin \phi)^{2\l}
           d\phi,
\end{align*}
where $q_\phi(\cos\ta)=q(\cos\phi,\cos \ta)$.  Since $\mu'\in
(0,1]$, $q_\phi(\cos \phi)=(-1)^r q_{\phi }(-\cos\phi) =1$ and
$\sup_{\phi, \ta} |q_\phi'(\cos\ta)|\leq c<\infty$, the desired
lower  estimate in this case follow by a slight modification of
the proofs in Cases 1 and 2.

\medskip

Putting these  cases together, we have completed the proof of
Proposition \ref{lowerJacobi}.

\enddocument